\documentclass[graybox]{svmult}
\usepackage{mathptmx}       
\usepackage{helvet}         
\usepackage{courier}        
\usepackage{makeidx}         
\usepackage{graphicx}        
\usepackage{multicol}        
\usepackage[bottom]{footmisc}
\usepackage{amsmath,amssymb,amsfonts}        
\makeindex             
\usepackage{amsmath,amssymb,amsfonts,epsfig,color,euscript}

\begin{document}


\newcommand{\ig}{\includegraphics}
\newcommand{\rb}{\raisebox}
\newcommand\risS[6]{\rb{#1pt}[#5pt][#6pt]{\begin{picture}(#4,15)(0,0)
  \put(0,0){\ig[width=#4pt]{#2.eps}} #3
     \end{picture}}}
\newcommand\risSpdf[6]{\rb{#1pt}[#5pt][#6pt]{\begin{picture}(#4,15)(0,0)
  \put(0,0){\ig[width=#4pt]{#2.pdf}} #3
     \end{picture}}}
\def\vyd#1{\mbox{\underline{\it\bfseries #1}}}
\def\kb#1{[ #1 ]}
\newcommand{\A}{\mathcal{A}}  
\def\a{\alpha}
\def\b{\beta}
\def\Ga{G}
\def\de{\delta}
\def\s{\sigma}
\def\S{\EuScript {S}}
\def\bc{\mathrm{bc}}
\def\I{\mathcal I}
\def\C{\mathcal C}
\def\B{\mathcal B}
\def\leq{\leqslant}
\def\geq{\geqslant}
\newcommand{\R}{\mathbb{R}} 
\newcommand{\Co}{\mathbb{C}} 
\def\red#1{\begin{color}{red}#1\end{color}}
\def\blue#1{\begin{color}{blue}#1\end{color}}
\def\yellow#1{\begin{color}[RGB]{255,255,0}{#1}\end{color}}
\def\magenta#1{\begin{color}[RGB]{255,0,255}{#1}\end{color}} 
\def\brown#1{\begin{color}[RGB]{192,96,0}{#1}\end{color}} 
\newcommand\cd[1]{\risS{-9}{#1}{}{24}{15}{15}}
\newtheorem{thm}{Theorem}[section]
\newtheorem{defn}[thm]{Definition}
\newtheorem{exa}[thm]{Example}
\newtheorem{exas}[thm]{Examples}
\newtheorem{prop}[thm]{Proposition}
\newtheorem{cor}[thm]{Corollary}
\newtheorem{subs}[thm]{}   


\title*{Partial duality for ribbon graphs}
\author{Sergei Chmutov}
\institute{
The Ohio State University at Mansfield,
1760 University Drive Mansfield, OH 44906, USA.
{\tt chmutov.1@osu.edu}}


\maketitle
\abstract{This is an expository paper extending the tutorial talk at the MATRIX Workshop on Uniqueness and Discernment in Graph Polynomials in October 2023. The explanation is mainly based on \cite{CVT22} with maximal possible simplifications.}

\section{Introduction} \label{s:intro}
The concept of partial duality appeared in \cite{Ch} in an attempt to unify the various Thistlethwaite type theorems of \cite{CP,DFKLS,CV} from the knot theory.
The classical Thistlethwaite theorem relates the Jones polynomial of a link to the Tutte polynomial of a checkerboard graph, the Tait graph, of its plane diagram. There were several generalizations of this theorem to virtual links where graphs on a higher genus surface, i.e. ribbon graphs, naturally appeared. It turned out that all these ribbon graphs are related via partial duality operations.

Since then partial duality was thoughtfully studied and described in many papers \cite{CVT21,EZ,GMT,GMT21,GMT21e,Mo3,YaJi21,YaJi22,YaJi22,YaZh22,YaJi22s,YaJi22dm,Yus22}. For an excellent exposition see \cite{EMM}. In \cite{CMNR19c,CMNR19p}, it was found that the partial dual operation generalized to delta matroids is a well known operation of twist of A.~Bouchet \cite{Bu87,Bu89}. So we should treat the partial duality from \cite{Ch} simply as a geometric description of Bouchet's twist operation for the delta matroids of ribbon graphs.
Also we would like to mention \cite{CVT22}, where the partial duality operation was generalized to hypermaps. Its (delta) matroidal analogue is still unknown.

The idea of partial duality is to split the
classical {\it Euler-Poincar\'e} duality
$$G\ =\ \risS{-20}{tr}{\put(3,20){\brown{$\scriptstyle e_1$}}
                 \put(25,20){\brown{$\scriptstyle e_2$}}
								 \put(20,-4){\brown{$\scriptstyle e_3$}}
              }{50}{10}{0}
\quad \risS{-5}{totor}{}{30}{20}{30}\quad
\risS{-30}{tr-g2}{}{60}{0}{0}
\quad \risS{-5}{totor}{}{30}{0}{0}\quad
\risS{-20}{g2}{\put(5,20){\magenta{$\scriptstyle e_1$}}
               \put(47,20){\magenta{$\scriptstyle e_2$}}
							 \put(20,12){\magenta{$\scriptstyle e_3$}}
              }{60}{0}{0}\quad =\ G^*
$$
for graphs cellulary embedded to surfaces into steps making a duality relative to a single edge at a time. Before giving the formal definitions in the next section we explain how it works on this particular example. The difficulty is that the partial duality may change the genus of the surface into which the graph is embedded.

First, to find the partial duality relative to the edge $e_1$ it is convenient to consider the graph $G$ embedded not into a plane but into an annulus, a surface with boundary.
$$\risS{-20}{tr-s}{\put(-5,35){$G$}\put(8,25){\brown{$\scriptstyle e_1$}}
                 \put(31,25){\brown{$\scriptstyle e_2$}}
								 \put(20,6){\brown{$\scriptstyle e_3$}}
              }{45}{20}{20}$$
The partial duality relative to the edge $e_1$ changes the annulus by adding a band across the edge $e_1$ to the two boundary circles of the annulus. Then the edge $e_1$ is shrunk to a new vertex and the dual edge $e_1$ goes along the band. Here is the result $G^{\{e_1\}}$, the partial dual of $G$ relative to the edge $e_1$.
$$\risS{-20}{tr-e1}{\put(50,55){$G^{\{e_1\}}$}
                 \put(18,40){\magenta{$\scriptstyle e_1$}}
                 \put(55,30){\brown{$\scriptstyle e_2$}}
								 \put(36,5){\brown{$\scriptstyle e_3$}}
              }{75}{45}{25}
$$
Notice that this is a duality. This means that if we would like to make a partial dual relative to a loop $e_1$ of the graph $G^{\{e_1\}}$ then we have to return to the original graph $G$. Hence, the partial duality relative to a loop removes the band along this loop and split the corresponding vertex into two vertices connected by the dual edge. The partial duality of $G^{\{e_1\}}$ relative to the edge $e_2$ works in a similar way to produce the graph $G^{\{e_1,e_2\}}$.
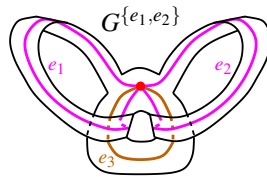
$$\risS{-20}{tr-e1e2}{\put(35,55){$G^{\{e_1,e_2\}}$}
                 \put(15,40){\magenta{$\scriptstyle e_1$}}
                 \put(78,40){\magenta{$\scriptstyle e_2$}}
								 \put(34,5){\brown{$\scriptstyle e_3$}}
              }{100}{45}{25}
$$
Now the edge $e_3$ is a loop. According to the above procedure for partial duality relative to a loop we have to cut the band along this loop and split the only vertex of $G^{\{e_1,e_2\}}$. The result is the graph
$G^{\{e_1,e_2,e_3\}}$
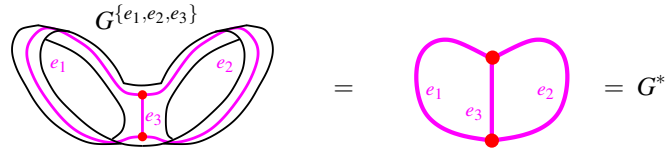
$$\risS{-20}{tr-e1e2e3}{\put(32,45){$G^{\{e_1,e_2,e_3\}}$}
                 \put(15,30){\magenta{$\scriptstyle e_1$}}
                 \put(78,30){\magenta{$\scriptstyle e_2$}}
								 \put(51,10){\magenta{$\scriptstyle e_3$}}
              }{100}{35}{25}\qquad=\qquad
\risS{-20}{g2}{\put(5,20){\magenta{$\scriptstyle e_1$}}
               \put(47,20){\magenta{$\scriptstyle e_2$}}
							 \put(20,12){\magenta{$\scriptstyle e_3$}}
              }{60}{0}{0}\quad =\ G^*
$$
on a sphere which is equivalent to the classical Euler-Poincar\'e dual $G^*$.

\section{Formal definitions} \label{s:form_def}
\begin{defn}\label{def:rg}{\rm
A {\it\bfseries ribbon graph} $G$ is a surface (possibly non-orientable) with boundary, represented as the union of two sets of closed topological discs called {\it\bfseries vertex-discs} $V(G)$ and 
{\it\bfseries edge-ribbons} $E(G)$, satisfying the following
conditions:
\begin{itemize}
\item[$\bullet$] these discs intersect by disjoint line segments;
\item[$\bullet$] each such line segment lies on the boundary of precisely one vertex-disc and precisely one edge-ribbon;
\item[$\bullet$] every edge-ribbon contains exactly two such line segments.
\end{itemize} }
\end{defn}

We consider ribbon graphs up to a homeomorphism of the corresponding surfaces preserving the decomposition on vertices and edges.

\begin{defn}\label{def:pd}{\rm
For a ribbon graph $G$ and a subset of the edge-ribbons $A\subseteq E(G)$, the {\it\bfseries partial dual}, $G^A$ of $G$ relative to $A$ is a ribbon graph constructed as follows. 
\begin{itemize}
\item The vertex-discs of $G^A$ are bounded by connected components of the boundary of the spanning subgraph of $G$ containing all the vertices of $G$ and only the edges from $A$. 
\item The edge-ribbons of
$E(G)\setminus A$ are attached to these new vertices exactly at the same places as in $G$. The edge-ribbons from $A$ become parts of the new vertex-discs. 
For $e\in A$ we take a copy of $e$, $e'$, and attach it to the new vertex-discs in the following way. The rectangle representing $e$ intersects with vertex-discs of $G$ by a pair of opposite sides. But it intersects the boundary of the spanning subgraph, that is the new vertex-discs, along the arcs of the other pair of its opposite sides. We attach $e'$ to these arcs by
this second pair. The copies of the first pair of sides in $e'$ become the arcs of the boundary of $G^A$.
\end{itemize}}
\end{defn}

Partial duality relative to a set of edges can be done step by step one edge at a time. 

The partial duality relative to one edge can be illustrated as follows.
$$G =\ 
  \risS{-8}{dus1}{\put(23,-3){\mbox{$e$}}}{50}{0}{0}\quad 
   \risS{-4}{totorl}{}{20}{0}{0}\quad 
  \risS{-24}{dus3}{\put(-8,13){\mbox{$e'$}}}{75}{25}{30}\ =\ 
  \risS{-25}{dus4}{\put(31,50){\mbox{$e'$}}}{50}{0}{0}\ =\ G^{\{e\}}
$$
Here the boxes with dashed arcs mean that there might be other edges attached to these vertices. The partial duality relative to an orientable loop goes in the opposite direction, from right to left.

For a non-orientable loop we have.
$$G =\ 
  \risS{-20}{dus-no1}{\put(31,40){\mbox{$e$}}}{50}{0}{0}\quad 
   \risS{-2}{totorl}{}{20}{25}{30}\quad 
  \risS{-20}{dus-no2}{\put(31,40){\mbox{$e'$}}}{50}{0}{0}\ =\ G^{\{e\}}
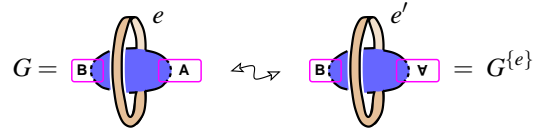
$$

\def\wt{\widetilde}
\begin{lemma}{\bf Properties of partial duality \cite{Ch}.} \label{lem:prop}
Let $G$ be a ribbon graph.
\begin{itemize}
\item[(a)] Suppose $E(G)\ni e\not\in A\subseteq E(G)$. Then 
$G^{A\cup\{e\}} = \bigl(G^A\bigr)^{\{e\}}$.
\item[(b)] $\bigl(G^A\bigr)^A= G$.
\item[(c)] $\bigl(G^A\bigr)^B= G^{\Delta(A,B)}$, where 
    $\Delta(A,B):=(A\cup B)\setminus (A\cap B)$ is the symmetric difference of sets.
\item[(d)] Partial duality preserves orientability of ribbon graphs.
\item[(e)] Let $\wt{G}$ be a surface without boundary obtained from $G$ by gluing discs to all boundary component of $G$. Then $\wt{G^A} = \wt{G^{E(G)\setminus A}}$.
\item[(f)] The partial duality preserves the number of connected components of ribbon graphs.
\end{itemize}
\end{lemma}

\bigskip
There are several ways to encode ribbon graphs. In subsequent sections we express partial duality relative to a single edge in these encodings following \cite{CVT22} where it was given in a more general situation of hypermaps.

\section{Partial duality for gems} \label{s:pd_gem}
\begin{defn}\label{def:gem}{\rm
According to S.~Lins \cite{Lins82} a {\it\bfseries gem (graph-encoded map)} is a trivalent graph whose edges are colored into three colors 0, 1, and 2 in such a way that that at every vertex the three meeting edges all have different colors.}
\end{defn}
Gems encode ribbon graphs in the following way. The trivalent graph is the graph of boundary segments of all the discs. The line segments between the vertex-discs and edge-ribbons are colored 2. The other segments of the boundary of vertex-discs (in between vertices and faces) colored 1. And the other two sides of edge-ribbons (in between edges and faces) are colored 0. Backwards, from a gem we can construct a ribbon graph in the following way. To each cycle of the gem colored by two colors 1 and 2, 12-cycle, we glue a disc which will be a vertex-disc of the resulting ribbon graph. To each 02-disc we glue an edge-ribbon. In general it will be a hypermap. But if all 02-cycles of a gem have length 4, then we will get a ribbon graph. This construction can be generalize to higher dimensions with more colors, see \cite{CVT22}.  

Here is an example.
$$G\ =\quad \risS{-15}{tr}{}{40}{20}{20}\quad
     =\quad \risS{-20}{tr-rg}{}{50}{0}{0}\quad
     =\quad \risS{-20}{tr-gem}{\put(40,35){\scriptsize gem}}{50}{0}{0}
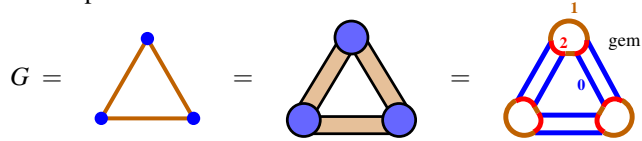
$$

\bigskip
\begin{prop}{\bf Partial duality relative to an edge for gems \cite{CVT22,EZ}.}\label{prop:pd_gem} 
Let $e$ be an edge of a ribbon graph $G$ and $C$ be the corresponding 
02-cycle of its gem.
Then the partial duality relative to $e$ is given by swapping colors 1 and 2 along the cycle $C$.
\end{prop}
Here is an example
$$G\ =\quad 
 \risS{-20}{tr-gem}{\put(24,-5){\brown{$\scriptstyle e_3$}}}{50}{25}{20}
\quad \risS{2}{totor}{}{20}{0}{0}\quad
 \risS{-20}{tr-gem-e3}{\put(24,-5){\magenta{$\scriptstyle e_3$}}}{50}{0}{0}
\quad =\quad
 \risS{-35}{tr-rg-e3}{}{55}{0}{0}\quad =\quad G^{\{e_3\}}
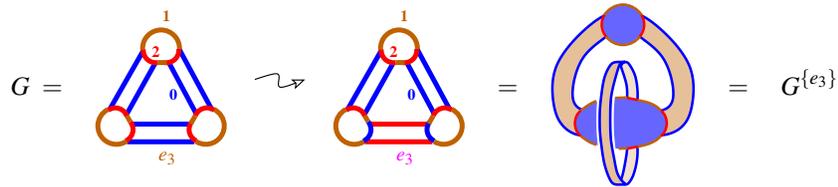
$$

\section{Partial duality for rotation systems} \label{s:pd_rs}

A rotation system is a way to encode a ribbon graph by two permutations of $2|E|$ elements.
\begin{defn}\label{def:rs} {\rm \cite{GT,MT}
A {\it\bfseries rotation system} on an abstract graph is cyclic order of half-edges around every vertex (given by a permutation $\sigma_V$) and an involution of two half-edges forming a single edge (given by a product of transpositions $\sigma_E$).}
\end{defn}
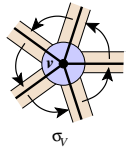
$$\risS{-20}{rot-sys}{}{50}{30}{10}
$$
Here is an example.
$$\rb{8pt}{$G=$}
\risS{-10}{tr}{\put(12,26){$\scriptstyle 1$}\put(1,9){$\scriptstyle 2$}
	               \put(6,-4){$\scriptstyle 3$}\put(28,-4){$\scriptstyle 4$}
								 \put(25,26){$\scriptstyle 6$}
								 \put(35,9){$\scriptstyle 5$}}{40}{0}{10}\qquad\qquad
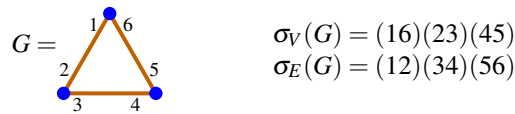
\begin{array}[b]{l}\sigma_V(G)=(16)(23)(45)\\ \sigma_E(G)=(12)(34)(56)
\end{array}
$$

We use rotation systems for orientable ribbon graphs. In that case, two permutations $(\sigma_V, \sigma_E)$ completely determine the ribbon graph.

\bigskip
\begin{prop}{\bf Partial duality relative to an edge for rotation systems \cite{CVT22}.}\label{prop:pd_rs} 
Suppose that an edge $e$ has two half edges $a$ and $b$. So, in the permutation $\sigma_E$ it corresponds to a transposition $(ab)$. Then
$$\sigma_V(G^{\{e\}}) = (ab)\cdot\sigma_V(G),\quad
\sigma_E(G^{\{e\}}) = \sigma_E(G)\ .$$
\end{prop}

These formulas were rediscovered in \cite[Equation 14]{GMT21e} and announced in \cite{GMT21} (see also \cite[Section 5.2]{GMT}).

For example, if  $e_3=(34)$ in our triangle example, then\\
$$\sigma_V(G^{\{e_3\}}) = (34)\cdot(16)(23)(45)=(16)(2453)$$.
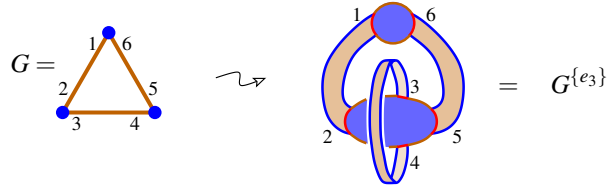
\vspace{-10pt}

$$\rb{8pt}{$G=$}
	\risS{-10}{tr}{\put(12,26){$\scriptstyle 1$}\put(1,9){$\scriptstyle 2$}
	               \put(6,-4){$\scriptstyle 3$}\put(28,-4){$\scriptstyle 4$}
								 \put(25,26){$\scriptstyle 6$}
								 \put(35,9){$\scriptstyle 5$}}{40}{0}{0}\qquad
\risS{2}{totor}{}{20}{40}{20}\qquad
\risS{-35}{tr-rg-e3}{\put(12,62){$\scriptstyle 1$}
                 \put(1,16){$\scriptstyle 2$}
	               \put(34,35){$\scriptstyle 3$}\put(34,6){$\scriptstyle 4$}
								 \put(40,62){$\scriptstyle 6$}
								 \put(50,16){$\scriptstyle 5$}}{55}{0}{0}\quad =\quad G^{\{e_3\}}
$$

\section{Partial duality for bi-rotation systems} \label{s:pd_rs}

A bi-rotation system is a way to encode a ribbon graph by three permutations $\tau_0$, $\tau_1$, and $\tau_2$ of $4|E|$. It works for non-orientable ribbon graphs as well.
\begin{defn}\label{def:brs} {\rm \cite{GT,MT}}
A {\it\bfseries bi-rotation system} is a set of three fixed point free involutions, $\tau_0$, $\tau_1$, and $\tau_2$,  acting on a set of vertices of the corresponding gem according to the figure.
$$\risS{-20}{flag}{}{70}{15}{15}\qquad
\risS{-20}{tau-0}{}{80}{0}{0}\qquad
\risS{-20}{tau-1}{}{60}{0}{0}\qquad
\risS{-20}{tau-2}{}{60}{0}{0}
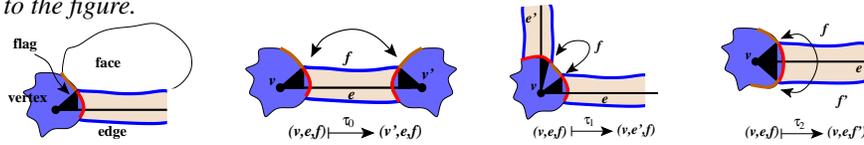
$$
\end{defn}
Here we think about vertices of the gem as {\it\bfseries local flags} consisting of a triple $(v,e,f)$ of a vertex $v$, an edge $e$ incident to $v$, and a face $f$ adjacent to both $v$ and $e$. We will depict such a flag as a small black right triangle with one acute angle at the vertex $v$, the right angle at the edge $e$ and the another acute angle at the vertex of the gem. So every edge-ribbon of the original ribbon graph has fours such flags adjacent to it.

At each vertex of the gem the three edges of the gem are colored by three different colors, $0$, $1$, and $2$. The permutation $\tau_i$ swaps the flags in the pairs of vertices connected by an an edge colored by $i$. 

In \cite{Tutte1984aa} W.~Tutte introduced a less symmetrical description of a ribbon graph (combinatorial map) in terms of three permutations $\theta$, $\phi$, and $P$ which can be expressed in terms of $\tau_0$, $\tau_1$, and $\tau_2$ as follows:
$$\theta=\tau_2,\quad \phi=\tau_0,\quad P=\tau_1\tau_2.
$$

Here is an example of the bi-rotation system.
$$G=\risS{-20}{tr-gem-bir}{\put(12,35){$\scriptstyle 1$}
    \put(1,16){$\scriptstyle 2$}\put(15,12){$\scriptstyle 3$}
		\put(21,34){$\scriptstyle 4$}
		\put(8,5){$\scriptstyle 5$}
    \put(10,-5){$\scriptstyle 6$}\put(35,-5){$\scriptstyle 7$}
		\put(32,10){$\scriptstyle 8$}
		\put(39,8){$\scriptstyle 9$}
    \put(45,15){$\scriptstyle a$}\put(33,37){$\scriptstyle b$}
		\put(25,25){$\scriptstyle c$}
		\put(23,3.5){$\scriptstyle e$}
       }{50}{0}{28}\qquad\qquad\qquad
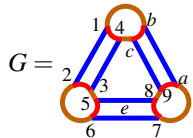
\begin{array}{l}
\tau_0(G)=(12)(34)(58)(67)(9c)(ab)\\
\tau_1(G)=(1b)(4c)(26)(35)(7a)(89)\\
\tau_2(G)=(1c)(23)(56)(78)(9a)(bc)
\end{array}
$$
In this case we use the hexadecimal number system. So, $a=10$, $b=11$, and $c=12$. In particular, for the edge-ribbon $e$ we have four local flags 5,6,7,and 8. So its contribution to the permutations 
$\tau_0$ and $\tau_2$ will be $\tau_0^e=(58)(67)$ and 
$\tau_2^e=(56)(78)$.

\bigskip
\begin{prop}{\bf Partial duality relative to an edge for bi-rotation systems \cite{CVT22}.}\label{prop:pd_brs} 
Let $\tau_0^e$ be the product of two transpositions in $\tau_1$ for $e$, and $\tau_2^e$ be the product of two transpositions in $\tau_2$ for $e$.

\medskip
Then\quad
$\tau_0(G^{\{e\}})=\tau_0(G)\tau_0^e\tau_2^e,\quad 
\tau_1(G^{\{e\}})=\tau_1(G),\quad
\tau_2(G^{\{e\}})=\tau_2(G)\tau_0^e\tau_2^e$.
\end{prop}

In our example we have
$$\begin{array}{ll}
\tau_0(G^{\{e\}})&=\tau_0(G)\tau_0^e\tau_2^e
  =(12)(34)(58)(67)(9c)(ab)(58)(67)(56)(78)\\ &=(12)(34)(56)(78)(9c)(ab)\\[5pt]
\tau_1(G^{\{e\}})&=(1b)(4c)(26)(35)(7a)(89)\\[5pt]
\tau_2(G^{\{e\}})&=\tau_2(G)\tau_0^e\tau_2^e
  =(1c)(23)(56)(78)(9a)(bc)(58)(67)(56)(78)\\ &=(1c)(23)(58)(67)(9a)(bc)
\end{array}
$$
These permutations correspond to the following gem:
$$	
G^{\{e\}}=\risS{-20}{tr-gem-e3-bir}{\put(12,35){$\scriptstyle 1$}
    \put(1,16){$\scriptstyle 2$}\put(15,12){$\scriptstyle 3$}
		\put(21,34){$\scriptstyle 4$}
		\put(8,5){$\scriptstyle 5$}
    \put(10,-5){$\scriptstyle 6$}\put(35,-5){$\scriptstyle 7$}
		\put(32,10){$\scriptstyle 8$}
		\put(39,8){$\scriptstyle 9$}
    \put(45,15){$\scriptstyle a$}\put(33,37){$\scriptstyle b$}
		\put(25,25){$\scriptstyle c$}
		\put(23,3.5){$\scriptstyle e$}
       }{50}{20}{20}
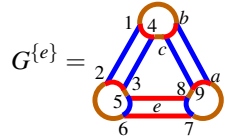
$$

\section{Genus change} \label{s:genus}

In this section we give a formula for the change of genus under partial duality.
We will deal with the {\it\bfseries Euler genus} $\gamma$ which is equal to twice the genus for orientable ribbon graphs and to the number of M\"obius bands $\mu$ in the presentation of the surface of the ribbon graph as (possibly several) spheres with 
$\mu$ cross-caps in the unorientable case. 

  Let $G$ be a ribbon graph and $A$ be a subset of its edges. Consider a subgraph of $G$, $G[A]$ induced by $A$. This means that the edge-ribbons of $G[A]$ are exactly $A$ and the vertex-discs are exactly the vertex-discs of $G$ incident with edges of $A$. In a similar way we can considered the subgraph $G^{*}[A]$ of the total Euler-Poincaré dual $G^{*}$ induced by $A$. Its vertices are the vertex-discs of $G^{*}$, which are the faces of the original ribbon graph $G$.
Note that the graphs $G[A]$ and $G^{*}[A]$ considered as ribbon graphs themselves are not dual to each other. The vertices of $G^{*}[A]$ are the faces of $G$, not 
$G[A]$. For any ribbon graph $G$, we use the notation $v(G)$ and $f(G)$ for the number of vertex-discs of $G$ and the number of face-discs of $G$ respectively. The last one is the same as the number of boundary components of $G$.

\bigskip
\begin{prop}\label{prop:pd_genus} 
$$\gamma(G^{A})-\gamma(G)=v(G[A])+v(G^{*}[A])-f(G[A])-f(G^{*}[A]).$$
\end{prop}

This proposition is a further simplification of \cite[Corollary 3.4]{CVT22}. It covers various cases of \cite[Table 1.1]{GMT} and \cite[Fig.1]{CVT21}.

In our example
$$G\ =\quad \risS{-20}{tr-rg}{\put(20,-5){\scriptsize $e_3$}}{50}{0}{20}\quad
$$
with $A=\{e_3\}$, $G[A]$ is the ribbon graph with two vertices and an edge $e_3$ between them. So $v(G[A])=2$, but $f(G[A])=1$ because $G[A]$ has only one boundary component. For $G^{*}[A]$, there are two vertices of $G^{*}$ bacause there two faces of $G$. And the edge $e_3$ connects these two vertices of $G^{*}[A]$. So we have $v(G^{*}[A])=2$ and $f(G^{*}[A])=1$. A reader can see that the graphs $G[A]$ and $G^{*}[A]$ are equivalent, but not dual to each other as ribbon graphs! Therefore the change of genus from Proposition \ref{prop:pd_genus} is
$$\gamma(G^{\{e_3\}})-\gamma(G)=2+2-1-1=2.$$
This is an orientable case, the genus is half of the Euler genus $\gamma$. Consequently the graph
$G^{\{e_3\}}$ is on torus as we saw before.

\section{Partial-dual Genus Polynomial} \label{s:pd-pol}
In \cite{GMT} J.~L.~Gross, T.~Mansour, T.~W.~Tucker introduced the Partial-dual Genus Polynomial of ribbon graphs. In the orientable case it is as follows.
$${\ }^\partial\Gamma_G(z):=\sum_{A\subseteq E(G)} z^{g(G^A)}\ .$$
In the non-orientable case you need to use the Euler genus instead of genus.
This polynomial was studied in various papers 
\cite{CVT21,GMT,GMT21,GMT21e,YaJi21,YaJi22,YaJi22,YaZh22,YaJi22s,YaJi22dm,Yus22}.
In particular in \cite{YaJi22dm,Yus22} it was generalized to delta-matroids. Among the most recent results I would like to mention \cite{Ch23}, where it was discovered that the Partial-dual Genus Polynomial of ribbon graphs appearing from chord diagrams
is a weight system from the theory of Vassiliev knot invariants. Recently (December 2023), Iain Moffatt informed me that he proved that Partial-dual Genus Polynomial weight system comes from the Lie algebra $\mathfrak{gl}_N$ weight system. Also Q.Deng, F.Dong, X.Jin, Q.Yan \cite{DDJY} generalize this result to the farmed chord diagrams.
\begin{acknowledgement}
I am very grateful to the MATRIX Institute for providing excellent research atmosphere for the Workshop on Uniqueness and Discernment in Graph Polynomials in October 2023 and to the MATRIX-Simons Travel Grant for covering a big part of my travel expenses.
\end{acknowledgement}


\bigskip
\begin{thebibliography}{ABCDE}
\bibitem{Bu87} A.~Bouchet, {\it Greedy algorithm and symmetric matroids}, Math. Program. {\bf 38} (1987) 147--159.

\bibitem{Bu89} A.~Bouchet, {\it A. Bouchet, Maps and delta-matroids}, Discrete Math. {\bf 78} (1989) 59--71.

\bibitem{Ch} S.~Chmutov,  {\it Generalized duality for graphs on surfaces and the signed Bollob\'as-Riordan polynomial}, 
   Journal of Combinatorial Theory, Ser. B {\bf 99}(3) (2009) 617--638.

\bibitem{Ch23} S.~Chmutov,  {\it Partial-dual genus polynomial as a weight system},
  Communications in Mathematics, {\bf 31}(3) (2023) 113--124.
	 doi.org/10.46298/cm.11709, 12 pages.

\bibitem{CP} S.~Chmutov, I.~Pak,  {\it The Kauffman bracket of virtual links and the      Bollob\'as-Riordan polynomial}, 
    Moscow Mathematical Journal, {\bf 7}(3) (2007) 409--418.
   
\bibitem{CVT21} S.~Chmutov, F.Vignes-Tourneret, 
   {\it On a conjecture of Gross, Mansour and Tucker}. 
   {\it European Journal of Combinatorics}, {\bf 97}(3) (2021) 
	 doi:j.ejc.2021.103368, 1--7.

\bibitem{CVT22} S.~Chmutov, F.Vignes-Tourneret, 
   {\it Partial Duality of Hypermaps}. 
   {\it Arnold Mathematical Journal}, (2022) doi:s40598-021-00194-8, 1--24.

\bibitem{CV} S.~Chmutov, J.~Voltz,  {\it Thistlethwaite's theorem for virtual links}, Journal of Knot Theory and its Ramifications, {\bf 17}(10) (2008) 1189--1198.

\bibitem{CMNR19c} C.~Chun, I.~Moffatt, S.~Noble, R.~Rueckriemen, {\it Matroids, Delta-matroids and Embedded Graphs}, Journal of Combinatorial Theory, Series A, {\bf 167} (2019) 7--59. (journal / arXiv)

\bibitem{CMNR19p} C.~Chun, I.~Moffatt, S.~Noble, R.~Rueckriemen, {\it On the interplay between embedded graphs and delta-matroids}, Proceedings of the London Mathematical Society, {\bf 118} (2019) 675--700. 

\bibitem{DFKLS} O.~Dasbach, D.~Futer, E.~Kalfagianni, X.-S.~Lin, N.~Stoltzfus, 
   {\it The Jones polynomial and dessins d'enfant}, Journal of Combinatorial Theory, Ser.B, {\bf 98} (2008) 384--399.

\bibitem{DDJY} Q.~Deng, F.~Dong, X.~Jin, Q.~Yan, {\it Partial-dual polynomial as a 
   framed weight system}. Preprint {\tt arXiv:2402.03799 [math.CO]}.

\bibitem{EZ} M.~Ellingham, X.~Zha, \it{Partial duality and closed 2-cell embeddings}, Journal of Combinatorics {\bf 8}(2) (2017) 227--254.

\bibitem{EMM} J.~A.~Ellis-Monaghan, I.~Moffatt, {\it Graphs on  
   Surfaces: Dualities, Polynomials, and Knots}, Springer, 2013.

\bibitem{GT} J.~L.~Gross and T.~W.~Tucker,
   {\it Topological graph theory}, Wiley, NY, 1987.

\bibitem{GMT} J.~L.~Gross, T.~Mansour, T.~W.~Tucker, {\it Partial 
  duality for ribbon graphs, I: Distributions}, 
	European Journal of Combinatorics {\bf 86} (2020) 103084, 1--20.

\bibitem{GMT21} J.~L.~Gross, T.~Mansour, T.~W.~Tucker, {\it Partial 
  duality for ribbon graphs, II: Partial-twuality polynomials and monodromy
  computations}, 
	European Journal of Combinatorics {\bf 95} (2021) 103329, 1--28.

\bibitem{GMT21e} J.~L.~Gross, T.~Mansour, T.~W.~Tucker, {\it  
  Enumerating Graph Embeddings and Partial-Duals by Genus and Euler Genus}, 
	Enumerative Combinatorics and Applications {\bf 1}(1) (2021) \#S2S1, 1--17.

\bibitem{Lins82} S.~Lins, {\it Graph-encoded maps}, Journal of Combinatorial Theory, Ser. B, {\bf 32}(2) (1982) 171--181.

\bibitem{Mo3} I.~Moffatt, {\it Separability and the genus of a partial dual}, 
    European Journal of Combinatorics {\bf 34} (2013) 355--378.
		
\bibitem{MT} B.~Mohar, C.~Thomassen, {\it Graphs on Surfaces},
   The Johns Hopkins University Press, 2001.

\bibitem{Tutte1984aa} W.~Tutte, {\it Graph Theory}, In: Rota, G.C. (ed.) {\it Encyclopedia of Mathematics and its Applications},
{\bf 21}. Addison-Wesley Publishing Company, Boston (1984).

\bibitem{YaJi21} Q.~Yan, X.~Jin, {\it Counterexamples to a 
  conjecture by Gross, Mansour, and Tucker on partial-dual genus  
	polynomials of ribbon graphs},
	European Journal of Combinatorics {\bf 93} (2021) 103285, 1--12.

\bibitem{YaJi22} Q.~Yan, X.~Jin, {\it Counterexamples to the 
  interpolating conjecture on partial-dual genus polynomials of ribbon
  graphs},
	European Journal of Combinatorics {\bf 102} (2022) 103493, 1--7.

\bibitem{YaJi22s} Q.~Yan, X.~Jin, {\it Partial-dual polynomials 
  and signed intersection graphs},
	Forum of Mathematics, Sigma {\bf 102}(e69)(2022) 1--16.

\bibitem{YaJi22dm} Q.~Yan, X.~Jin, {\it Twist monomials of binary  
  delta-matroids}.
  Preprint \verb#arXiv:2205.03487v1 [math.CO]#.

\bibitem{YaZh22} Y.~Yang,  X.~Zha, {\it Partial-dual Euler-genus 
  distributions for bouquets with small Euler genus},
	Ars Mathematica Contemporanea. {\bf 22}(3) (2022) \#P3.09, 13pp.

\bibitem{Yus22} D.~Yuschak, {\it Delta-matroids with twist monomials}.
  To appear in European Journal of Combinatorics. 
	Preprint \verb#arXiv:2208.13258v1 [math.CO]#. 
	


\end{thebibliography}
\end{document}